\documentclass[11pt]{article}
\usepackage[active]{srcltx}
\usepackage{epsfig,amssymb}

\begin{document}

\title{
Flat 3-webs via semi-simple Frobenius 3-folds}

\author{{\Large Sergey I. Agafonov} \\
\\
Department of Mathematics
\\
Federal University of Paraiba\\Jo\~ao Pessoa, Brazil\\
e-mail: {\tt sergey.agafonov@mat.ufpb.br} }
\date{}
\maketitle

\unitlength=1mm

\newtheorem{theorem}{Theorem}
\newtheorem{proposition}{Proposition}
\newtheorem{lemma}{Lemma}
\newtheorem{corollary}{Corollary}
\newtheorem{definition}{Definition}
\newtheorem{example}{Example}

\pagestyle{plain}

\begin{abstract}
\noindent

We construct flat 3-webs via semi-simple geometric Frobenius manifolds of dimension three
 and give geometric interpretation of the Chern connection of the web. These webs turned out to be biholomorphic to the characteristic webs on the solutions of the corresponding associativity equation. We show that such webs are hexagonal and admit at least one infinitesimal symmetry at each singular point. Singularities of the web are also discussed.
(secondary).
\end{abstract}

\section{Introduction}
B.Dubrovin introduced the notion of Frobenius manifold  geometrizing the physical context of WDVV-equations that arise originally in the two-dimensional topological field theory (see \cite{Df} and \cite{Mf}).
Frobenius manifold is a complex analytic manifold $M$ equipped with the following
analytic objects:

\begin{itemize}
\item commutative and associative multiplication on $T_pM$,
\item invariant non-degenerate flat inner product : $<u\cdot v,w>=<u,v\cdot w>$,
\item constant unity vector field $e$: $\nabla e=0$, $e \cdot v=v\ \forall v\in TM$,
\item linear Euler vector field  $E$: $\nabla (\nabla E)=0$,
\end{itemize}
satisfying the following conditions:
\begin{enumerate}
\item one-parameter group generated by $E$ re-scales the
multiplication and the inner product,
\item 4-tensor $(\nabla _z c)(u,v,w)$ is symmetric in
$u,v,w,z$, where $$c(u,v,w):=<u\cdot v,w>.$$
\end{enumerate}
In the above definition the symbol $\nabla$ stands for the Levi-Civita connection of the inner product $<\ ,\ >$.

\begin{figure}[th]
\hspace{-0.5cm} \epsfig{file=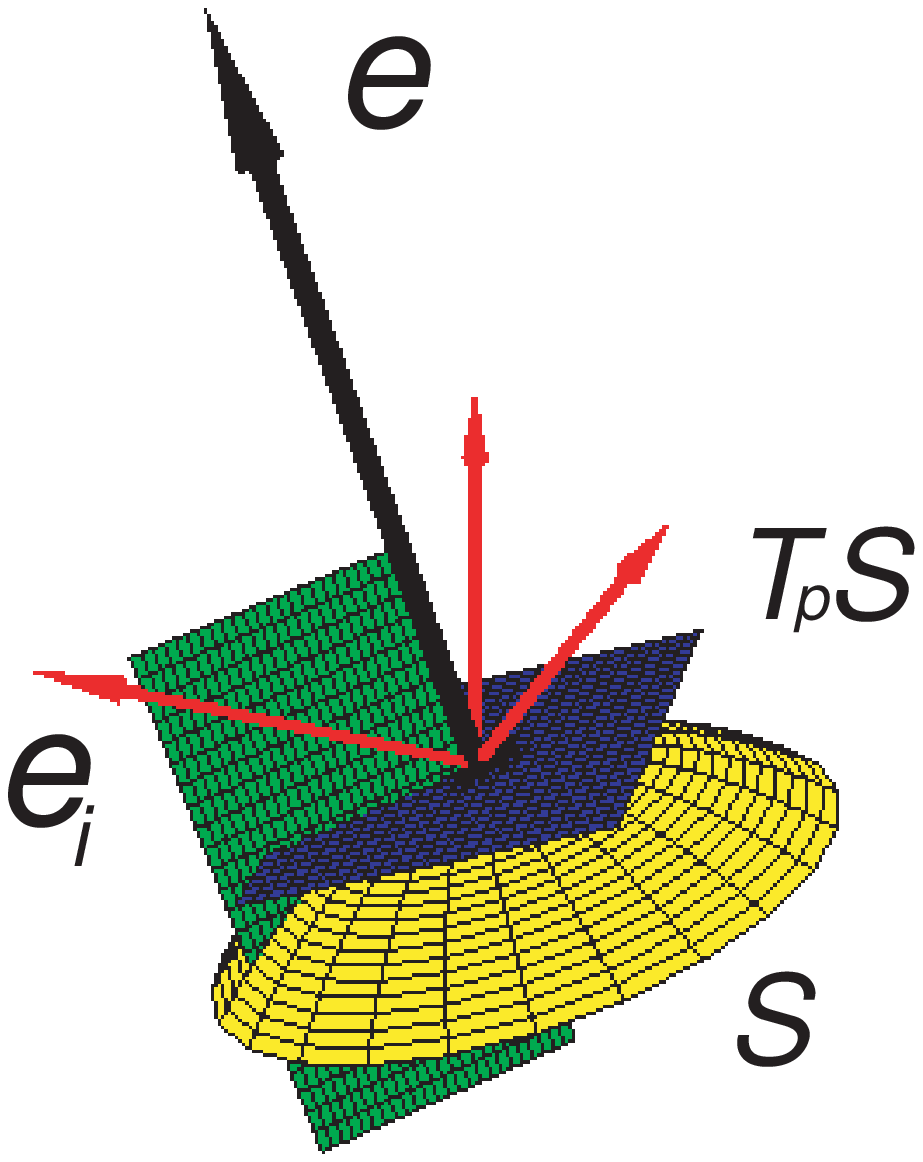,width=80mm} \ \ \ \put(1,8){ \epsfig{file=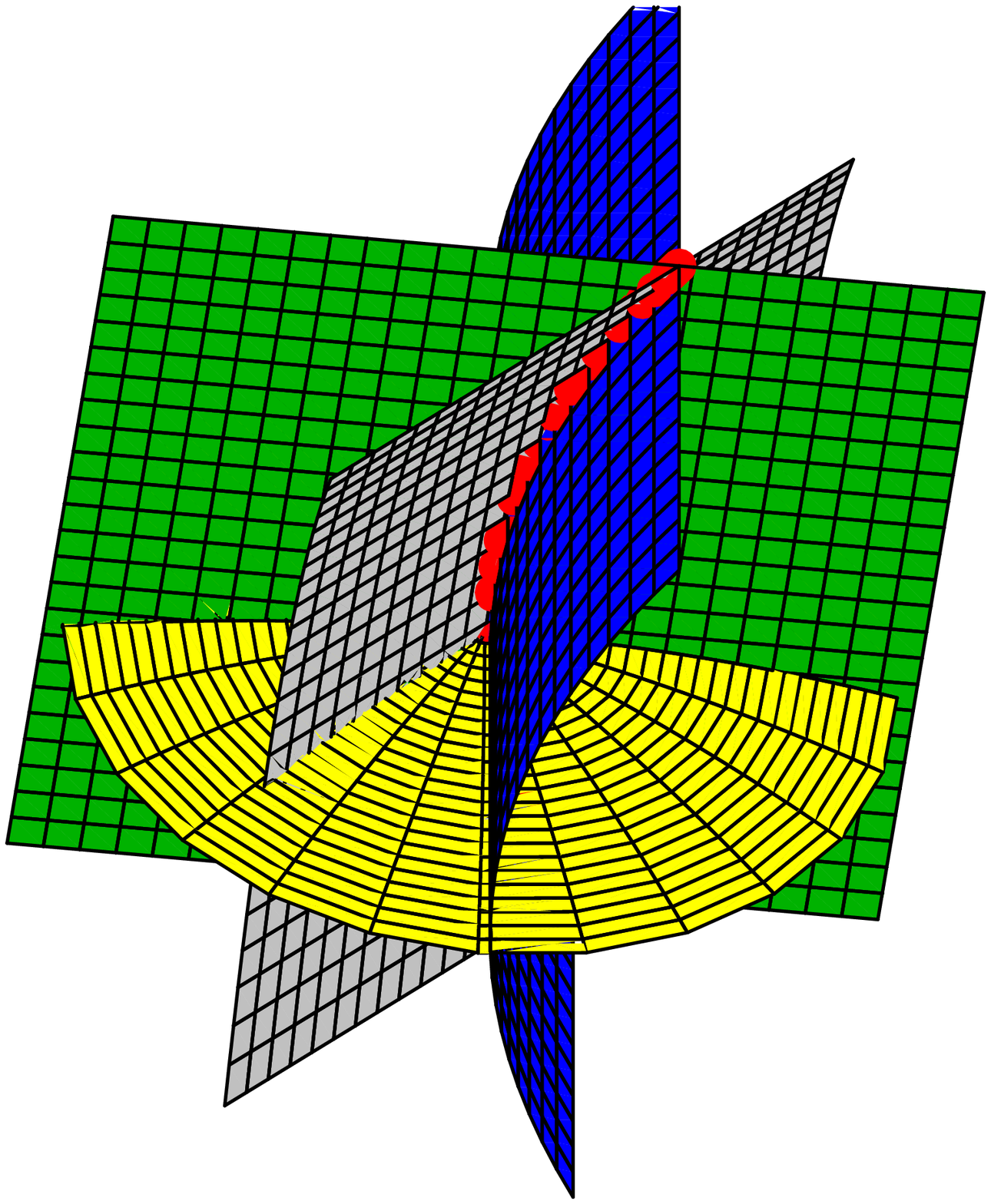,width=40mm}}  \caption{Construction of a booklet 3-webs from Frobenius 3-folds.} \label{booklet}
\end{figure}

Consider a semi-simple Frobenius manifold $M$ of dimension 3, which means that the algebra $T_pM$ is semi-simple for each $p\in U$ for some open set $U\subset M$. Then $T_pM$ is a direct
product of one-dimensional algebras spanned by idempotents:
$$T_pM=\mathbb C \{e_1\}\otimes \mathbb C \{e_2\}\otimes \mathbb C\{e_3\}, \ \ \ \ e_i\cdot e_j=\delta _{ij}.$$
In this setting the unity vector field is $$e=e_1+e_2+e_3.$$
Let $S$ be a surface transverse to the unit vector field $e$, then 3
planes spanned by $\{e,e_i\}$ cut 3 directions  on $T_pS$ (see Fig.\ref{booklet} on the left).
Integral curves of these direction fields build a {\bf flat or hexagonal 3-web}, i.e. this web is locally biholomorphic to 3 families of parallel lines.
Futher we will see that the distributions $\{e,e_i\}$ are integrable, integral surfaces of each of the distributions being formed by the integral curves of the unity vector field $e$. Thus for any point in $M$ there are 3 integral surfaces intersecting along a such curve.  These surfaces cut $S$ along the constructed web.  This justifies the following definition (see also Fig.\ref{booklet} on the right).
\begin{definition}
The constructed web is called a booklet 3-web.
\end{definition}
We call a point  {\bf regular or non-singular} if  the web directions are pairwise transverse.
Hexagonal 3-web does not have any local invariants in regular points. Its "personality"  is encoded in the behavior at {\bf singular points}, where at least 2 web directions coincide.

For the booklet 3-webs constructed above we show that:
\begin{itemize}
\item the booklet 3-web is biholomorphic to the characteristic 3-web of the corresponding solution of associativity equation,
\item the booklet 3-web has at least one infinitesimal symmetry at each singular point,
\item the Chern connection of the booklet 3-web is induced by the connection on $TM$ compatible with the algebraic structure of $M$.
\end{itemize}
Finally, we discuss possible singularities of the constructed webs.

\section{Characteristic Webs} Consider in more details the local construction described in the introduction. Suppose the operator of multiplication by the Euler vector field
$$
\mu_{E}: T_pM\to T_pM, \ \ v\mapsto E\cdot v
$$
has pairwise distinct eigenvalues $\lambda_1(p),\lambda_2(p),\lambda_3(p)$. Then the functions $\lambda_i(p)$ are local  coordinates called also {\bf canonical} or {\bf Dubrovin} coordinates.    In these coordinates $\partial_i \cdot \partial_j=\delta _{ij}\partial_i$, where $\partial_i =\partial / \partial _{\lambda_i}$ (see \cite{Di}). Therefore the basis vectors $\partial_i$ are idempotents $\partial_i=e_i$. Moreover, holds true $E=\sum_i\lambda_i\partial_i$.  The distributions spanned by the pairs $\{e_i,e\}$, where $e$ is the unity $e=e_1+e_2+e_3$, are determined by the following 1-forms:
\begin{equation}\label{thetas}
\theta_1=d\lambda_2-d\lambda_3, \ \ \theta_2=d\lambda_3-d\lambda_1, \ \ \theta_3=d\lambda_1-d\lambda_2.
\end{equation}
It is immediate that $\sum_i\theta_i=0$ and $d\theta_i=0$.
Let $\imath:U\in \mathbb C^2 \to S\subset M$ be a local parametrization of a surface $S$ transverse to the unity vector field $e$. Define $\omega_i:=\imath^*(\theta _i)$, then
$$
\sum_i\omega_i=\sum_i\imath^*(\theta_i)=\imath^*\left(\sum_i\theta_i\right)=0, \ \ d\omega_i=d(\imath^*(\theta _i))=\imath^*(d(\theta _i))=0.
$$
Thus the following Theorem is immediate.
\begin{theorem}
The booklet 3-web is hexagonal.
\end{theorem}
{\it Proof:} Each of the distributions $\omega_i=0$ on $S$ defines locally a first integral $u_i$ by $du_i=\omega_i$. These integrals sum up to a constant: $u_1+u_2+u_3=const$, which implies hexagonality of the web. (Take $u_1$ and $u_2$ as local coordinates on $S$.) \hfill $\Box$

\smallskip

Suppose the Lie derivative of the unity vector field $e$ does not vanish: $\mathcal{L}_E(e)\ne 0$. Dubrovin \cite{Df} showed that  any Frobenius
structure with this property is locally defined by a potential $F(t), \ t=(t^1,...,t^n)$ such
that its 3d order derivatives $c_{\alpha \beta
\gamma}(t):=\frac{\partial ^3F(t)}{{\partial t^{\alpha}
\partial t^{\beta} \partial t^{\gamma}}}$ satisfy the following conditions.
\begin{enumerate}
\item Normalization: the matrix $\eta _{\alpha \beta }=c_{1 \alpha \beta
}(t)$ is constant and not degenerate. This matrix gives inner  product $<\partial _\alpha ,\partial _\beta>=\eta
_{\alpha \beta}$, the parameters $(t^1,...,t^n)$ being its flat coordinates.
\item Associativity: the functions $c_{\alpha \beta}
^{\gamma}:=\sum _i \eta ^{\gamma i }c_{i \alpha \beta }(t)$ define
a structure of associative algebra on $T_pM$:
$$
\partial _{t^{\alpha} }\cdot \partial
_{t^{\beta}  }=\sum_i c_{\alpha \beta} ^{i}\partial _{t^i},
$$
with $(\eta ^{\alpha  \beta} ):=(\eta _{\alpha  \beta})^{-1}.$
\item Homogeneity:
$F(c^{w_1}t^1,...,c^{w_n}t^n)=c^{w_F}F(t^1,...,t^n)$.
\end{enumerate}
The algebra is automatically commutative, its unity is $e=\partial _1$,
the Euler vector field is given by $E=\sum_i w_it^i\partial _i$.
Associativity condition manifests itself as a system
of nonlinear PDEs, called associativity or
Witten-Dijkgraaf-Verlinde-Verlinde equations.

For Frobenius 3-folds there is only one associativity equation of order three.
There are two non-equivalent normal forms.
\begin{example} {\rm If $<e,e>= 0$ then in some flat coordinates the matrix
$\eta$ is anti-diagonal
$$\eta=\left(
\begin{array}{ccc}
0&0&1\\
0&1&0\\
1&0&0
\end{array}
\right)
$$
and the potential has the form
$$
F(t)=\frac{1}{2}t^2y+\frac{1}{2}tx^2+f(x,y), {\ \rm where \ } \ \
t^1=t,\ t^2=x, \ t^3=y.
$$
 Associativity equation reads as
$$
f_{yyy}=f_{xxy}^2-f_{xxx}f_{xyy}.
$$
Recall that characteristic curves $\varphi (x_1,x_2)=const$ for a solution of a PDE $H(x_1,x_2,f,f_{1},...,f_J)=0$, where $J=(j_1,j_2)$ with $j_1+j_2\le 3$ is a multi-index and $f_J=\frac{\partial ^{|J|}f}{\partial x_1^{j_1} \partial x_2^{j_2}}$, are defined by the following implicit PDE
\begin{equation}\label{equation PDE characteristics}
\sum _{|J|=3}\frac{\partial H}{\partial f_J}\nabla \varphi ^J=0, \ \ \ \ \nabla\varphi ^J:=\varphi_1^{j_1}\varphi_2^{j_2}.
\end{equation}
For the above equation this PDE takes the form
$$
\varphi _y^3+f_{xxx}\varphi _y ^2\varphi_x-2f_{xxy}\varphi _y\varphi _x^2+f_{xyy}\varphi _x^3=0
$$
As the associativity equation is non-linear the web of characteristics (we will call it also a {\bf characteristic 3-web}) depends on the solution.
Solving this equation is equivalent to integration of an implicit cubic ODE.
Substituting $[\varphi _y:\varphi _x]=[dx:-dy]$ one arrives at
the following cubic {\bf binary form}:
\begin{equation}\label{binaryPoly}
dx ^3-f_{xxx}dx ^2dy-2f_{xxy}dxdy ^2 -f_{xyy}dy ^3=0.
\end{equation}
The multiplication table for the corresponding Frobenius algebra is
$$
\begin{array}{l}
  \partial _x \cdot \partial _x=f_{xxy}\partial _t+f_{xxx}\partial _x+\partial _y,\\
  \partial _x \cdot \partial _y =f_{xyy}\partial _t+f_{xxy}\partial _x,\\
  \partial _y \cdot \partial _y =f_{yyy}\partial _t+f_{xyy}\partial _x.
\end{array}
$$
}\end{example}
\begin{example} {\rm If $<e,e> \ne 0$ then in some flat coordinates the matrix
$\eta$ is
$$\eta=\left(
\begin{array}{ccc}
1&0&0\\
0&0&1\\
0&1&0
\end{array}
\right)
$$
and the potential has the form
$$
F(t)=\frac{1}{6}t^3+txy+f(x,y), {\ \rm where \ } \ \
t^1=t,\ t^2=x, \ t^3=y.
$$
 Associativity equation reads as
$$
f_{xxx}f_{yyy}-f_{xxy}f_{xyy}=1.
$$
The characteristic  3-web in $(x,y)$-plane is
defined by the PDE
$$
f_{xxx}\varphi _y^3-f_{xxy}\varphi _y ^2\varphi_x-f_{xyy}\varphi _y\varphi _x^2+f_{yyy}\varphi _x^3=0
$$
or by the cubic binary form:
\begin{equation}\label{binaryDet}
f_{xxx}dx ^3+f_{xxy}dx ^2dy-f_{xyy}dxdy ^2 -f_{yyy}dy ^3=0.
\end{equation}
The multiplication table is
$$
\begin{array}{l}
  \partial _x \cdot \partial _x=f_{xxy}\partial _x+f_{xxx}\partial _y,\\
  \partial _x \cdot \partial _y =\partial _t+f_{xyy}\partial _x+f_{xxy}\partial _y,\\
  \partial _y \cdot \partial _y =f_{yyy}\partial _x+f_{xyy}\partial _y.
\end{array}
$$
}
\end{example}

\smallskip

\noindent {\bf Remark 1.} Note that the coefficients of the binary forms (\ref{binaryPoly}) and (\ref{binaryDet}) never vanish simultaneously. For this reason, the web directions are well-defined everywhere. They can coincide thus giving multiple roots of the binary equations.

\smallskip

\noindent {\bf Remark 2.}  The above two associativity equations are equivalent with respect to some non-local coordinate transform, depending on solutions (see \cite{Al}). This transform preserves the hexagonality of the web.

\bigskip

\noindent It was observed (see \cite{MFa}) that the characteristic 3-webs are hexagonal. Now we show that they are diffeomorphic to the above constructed booklet webs.
\begin{theorem}
The booklet 3-web is biholomorphic to the characteristic 3-web of the corresponding solution of associativity equation.
\end{theorem}
{\it Proof:} Let us check this claim for the plane $t=const$. Set $N=T\partial _t+X\partial _x+Y\partial _y$ for an idempotent. The coefficients by $\partial _t,\partial _x,\partial _y$ of the vector equation $N\cdot N-N=0$ give 3 scalar equations. The coefficient by $\partial _y$ gives $T$. Substituting this expression into the coefficient by $\partial _x$ one gets the binary differential forms for $[dx:dy]=[X,Y]$ that coincide with (\ref{binaryPoly}) or (\ref{binaryDet}) respectively. For the general case choose a plane $P=\{(t,x,y):t=const\}$ and consider the flow generated by the unity vector field. This flow preserves the distributions $\theta_i=0$. For each point $p\in S$ of the surface the orbit of the point  cuts the plane $P$ at a point $\psi (p)\in P$, thus giving a desired biholomorphism $\psi$. \hfill $\Box$

\section{Infinitesimal Symmetries}

We say that a 3-web on a surface with local coordinates $x,y$ has an infinitesimal symmetry $$
X=\xi(x,y)\partial _x+\eta (x,y)\partial _y,
$$
if the local flow  of the vector field $X$ maps the web leaves to the web leaves.
Infinitesimal symmetries form a Lie algebra with respect to the
Lie bracket. Cartan  proved  (see \cite{Cg}) that at a regular point
a 3-web either does not have infinitesimal symmetries
(generic case), or has a one-dimensional symmetry algebra (then in
suitable coordinates it can be defined by the form
$$dx\cdot dy \cdot (dy+g(x+y)dx)=0, \ \ \ g\ne const$$ with the symmetry $\partial _y
-\partial _x$), or has a three-dimensional symmetry algebra (then
it is equivalent to the web defined by the form  $$dx\cdot dy \cdot (dy+dx)=0$$ with the
symmetry algebra generated by  $\{\partial _x,\
\partial _y, x\partial_x+y\partial_y\}$).
In the last case, when the symmetry algebra has the largest possible  dimension 3,
the 3-web  is  hexagonal. Note that not all symmetries survive at a singular point. The condition to have at least one-dimensional
symmetry at a singular point is not trivial. The following
binary equation $$dy^3-2x^2y(1+x^2)dydx^2+8x^3y^2dx^3=0$$ defines a flat 3-web but
looses all its symmetries at the singular point $(0,0)$.

\begin{theorem}
The booklet 3-web has at least one infinitesimal symmetry at each point.
\end{theorem}
{\it Proof:} The flow $\exp(aE)$ of the Euler vector field $E$ respects the distributions $\theta_i=0$ and induces an action $T_a$ on the surface $S$. This action preserves the web. Indeed, let $p\in S$ be a point of the surface and $C_p$ the orbit of $p$ under the  flow of the unity vector field $e$.
 Then the  action $T_a$ maps $p$ to the point of intersection of $S$ with the image of $C_p$ under $\exp(aE)$, i.e.  $T_a(p):=\exp(aE)C_p\cap S$.  \hfill $\Box$

\section{Chern Connection}

For computing Chern connection of 3-webs we use Blaschke's approach based on differential forms \cite{Be}.
 Let us introduce {\bf binary vector field} as a section of the vector bundle with the base manifold $U\in \mathbb C^2$ and the fiber of symmetric forms $V$ on $T^*_p U$, i.e. $V$ is a n-linear symmetric map  $(T^*_p U)^n \to \mathbb C$.  Then
equation (\ref{equation PDE characteristics}) has the form $V(d\varphi)=0$, where
\begin{equation}\label{binaryV}
V=a(x,y)(\partial _x)^3+b(x,y)(\partial _x)^2\partial _y+c(x,y)\partial _x(\partial _y)^2+r(x,y)(\partial _y)^3
\end{equation}
is a cubic binary vector field. (Do not confuse monomials in $\partial _x,\partial _y$ with the product in Frobenius algebra, defined in introduction!)

At a non-singular point the field $V$ can be factorized:
\begin{equation}\label{factorsV}
V=V_1V_2V_3, \ \ {\rm where} \ \ V_i=q_i(x,y)\partial _x-p_i(x,y)\partial _y, \ \ \ i=1,2,3
\end{equation}
and
$$
\sigma_i=p_idx+q_idy, \ \ \ i=1,2,3,
$$
are the roots of $V$, which means $V(\sigma_i)=0$.
Due to the homogeneity of $V$ these roots are defined up to the multiplication by a non-vanishing factor.
They can be normalized to satisfy the condition
\begin{equation}\label{normalization}
\sigma_1+\sigma_2 + \sigma _3=0.
\end{equation}
Upon introducing an "area" form by  $$\Omega =\sigma_1\wedge
\sigma_2=\sigma_2\wedge \sigma_3=\sigma_3\wedge
\sigma_1=(p_1q_2-p_2q_1)dx\wedge dy,$$
define the Chern
connection form  as
$$
\gamma:
=h_2\sigma_1-h_1\sigma_2=h_3\sigma_2-h_2\sigma_3=h_1\sigma_3-h_3\sigma_1,
$$
where $h_i$ are determined by $$d\sigma_i=h_i\Omega.$$
The web is flat iff  the connection form
is closed: $d(\gamma)=0$. This implies $d\sigma_i=\gamma \wedge \sigma_i$. Putting
\begin{equation}\label{dk}
dk=-\gamma k,
\end{equation}  we introduce first integrals $u_i$ of the foliations  (at least locally) at a
regular point  by
\begin{equation}\label{du}
du_1=k\sigma_1,\ \ du_2=k\sigma_2,\ \ du_3=k\sigma_3,
\end{equation}
satisfying the equation $u_1+u_2+u_3=0$, which implies hexagonality.

\smallskip

\noindent {\bf Remark.}  Let $\eta_1, \eta_2, \eta_3$ be germs of differential forms in $(\mathbb C^2,q_0)$ defining a flat 3 web and satisfying the following conditions:
\begin{itemize}
\item the forms are closed: $d(\eta_i)=0, \ \ i=1,2,3,$
\item the forms define the web: $\eta_i\wedge \sigma_i=0, \ \ i=1,2,3,$
\item the forms sum up to zero: $\eta_1 + \eta_2 + \eta_3=0,$
\end{itemize}
then these forms are proportional to $k\sigma_i$:  $\eta_i=c k\sigma_i, \ \ i=1,2,3$, $c=const$. One says that the space of abelian relations is one-dimensional for a hexagonal 3-web. I other words the first integrals summing up to zero are defined up to a constant factor.

\smallskip

\begin{proposition}
In the normalization (\ref{normalization}) the 3-web defined by the binary 3-vector (\ref{binaryV}) has the Chern connection
$$
\gamma =\frac{\gamma _1dx+\gamma_2dy}{3D},
$$
where
$$
D=18\, a\, b\, c\, r-27\,a^2\,r^2-4\,a\,c^3+b^2\,c^2-4\,b^3\,r
$$ is the discriminant of the equation $V(\sigma)=0$ and
$$
\begin{array}{l}
\gamma_1= \left(15\,b\,c\,r -27\,a\,r^2-4\,c^3 \right)a_x+6\,r \left( 3\,b\,r-c^2\right)a_y  \\
\ \ \ \   +2\,b \left(c^2- 3\,b\,r\right)  b_x +3\,r \left( b\,c-9\,a\,r \right) b_y\\
\ \ \ \  +b \left( 9\,a\,r -b\,c \right)  c_x+6\,r \left( 3\,a\,c-b^2 \right)  c_y\\
\ \ \ \ +2\,b \left( b^2 - 3\,a\,c\right) r_x+ \left( 9\,a\,b\,r-12\,a\,c^2+3\,b^2\,c \right) r_y
 \end{array}
$$
$$
\begin{array}{l}
\gamma_2= \left( 9\,a\,c\,r-12\,b^2\,r+3\,b\,c^2\right) a_x+2\,c \left(c^2- 3\,b\,r\right) a_y+\\
\ \ \ \  6\,a \left( 3\,b\,r-c^2 \right) b_x+ c \left( 9\,a\,r -b\,c \right) b_y+ \\
\ \ \ \  3\,a \left(9\,a\,r -b\,c \right) c_x+2\,c \left(b^2- 3\,a\,c \right) c_y\\
\ \ \ \  6\,a \left( 3\,a\,c-b^2\right) r_x+ \left( 15\,a\,b\,c-27\,a^2\,r-4\,b^3 \right) r_y.
 \end{array}
$$
\end{proposition}
{\it Proof:} With $p_3=-p_1-p_2$, $q_3=-q_1-q_2$ equation (\ref{factorsV}) is equivalent to the following 4 equations:
\begin{equation}\label{pqV}
\begin{array}{l}
q_1q_2(q_1+q_2)=-a,\ \ \  p_1p_2(p_1+p_2)=r, \\
p_1(q_2^2+2q_1q_2)+p_2(q_1^2+2q_1q_2)=b,\\
q_1(p_2^2+2p_1p_2)+q_2(p_1^2+2p_1p_2)=-c.
 \end{array}
\end{equation}
Differentiating each of these equations with respect to $x$ and $y$ one gets all the derivatives of $p_i,q_i$. Substitution of these expressions into the formula for $\gamma$ and simplification with the help of equations (\ref{pqV}) yields the resulting formula for $\gamma$. \hfill $\Box$
\begin{corollary}
In the normalization (\ref{normalization}) the characteristic 3-web of associativity equations is
\begin{equation}\label{associativity connection}
\gamma =-\frac{1}{6}d(\ln D),
\end{equation}
where the components  $a,b,c,r$  of the binary vector are given by equation (\ref{equation PDE characteristics}).
\end{corollary}
Using the form $\gamma$ Blaschke (see \cite{Be}) introduced a parallel transport of the tangent vectors. The construction is the following: take vector fields $v_i, \ \ i=1,2,3$ tangent to the web curves, normalize them to satisfy $v_1+v_2+v_3=0$, and choose, for example, $\{v_1,v_2\}$ as a frame on the tangent bundle. Then each vector field $\xi$ can be represented locally as $\xi=\xi^1v_1+\xi^2v_2$ and the corresponding connection is defined by the following forms on the tangent bundle:
$$
\Gamma^1=d\xi ^1-\gamma \xi^1, \ \ \ \ \Gamma^2=d\xi ^2-\gamma \xi^2,
$$
i.e. the vector field $\xi$ is parallel along a curve iff the forms $\Gamma ^i$ vanish along it.
In particular, if our 3-web is hexagonal then the vectors $v_i$ can be normalized to commute. In this normalization holds true $\gamma =0$ and the vector fields with $\xi ^i=const$ are parallel along any curve.

Now we give a geometric interpretation of the Chern connection of the characteristic web.  At a regular point the idempotents $\{e_1,e_2,e_3\}$ form a frame of the tangent bundle $TM$: $\eta =\eta ^1e_1+\eta ^2e_2+\eta ^3e_3$ for each vector field $\eta$ on $M$. The projections of two of them, say, of $e_1,e_2$ along the unity $e$  form a frame $\{v_1,v_2\}$ on $S$.
\begin{theorem}\label{connectionANDalgebra}
Suppose $v=\eta ^1e_1+\eta ^2e_2+\eta ^3e_3 \in T_pS$ and a curve $\alpha :I\mapsto S$ with  $\alpha (I)\in S, \ \ \alpha (0)=p$ does not passes through singular points of the booklet 3-web.  Consider the vector field $\eta\in TM$ along $\alpha$ such that the coordinates $\eta ^i$ are kept constant. The projection of $\eta$ into $T_{\alpha(t)}S$ along $e$ is the parallel transport defined by the Chern connection.
\end{theorem}
{\it Proof:} It is sufficient to check that the vector fields $v_1,v_2\in TS$ commute. Holds true $\theta _1(e_1)=0$, $\theta_2(e_1)=-1$, where the forms $\theta _i$ are defined by (\ref{thetas}). Setting $e_1=v_1+k_1e$ one obtains $\theta _1(e_1)=\theta _1(v_1)=\omega _1(v_1)=0$ and $\theta _2(e_1)=\theta_2(v_1)=\omega _2(v_1)=-1$. Similar $\omega _1(v_2)=1$ and $\omega _2(v_2)=0$. The forms $\omega_i$ are closed therefore the vector fields $v_1,v_2$ commute.   \hfill $\Box$

\bigskip

\noindent {\bf Remark.}  The connection on $TM$ defined by Theorem \ref{connectionANDalgebra} preserves the idempotents and therefore respects the algebraic structure of the Frobenius manifold.

\section{Singularities of Characteristic 3-webs}
In \cite{Ai} hexagonal 3-webs were studied via implicit cubic ODEs, obtained from a binary differential form
\begin{equation}\label{binary}
K_3(x,y)dy^3+K_2(x,y)dy^2dx+K_1(x,y)dydx^2+K_0(x,y)dx^3=0,
\end{equation}
where the coefficients $K_i$ do not vanish simultaneously, by deviding by $dx^3$ or $dy^3$. Namely, if the equation is brought to the form
without the
quadratic term
$$
p^3+A(x,y)p+B(x,y)=0,
$$
then its Chern connection is
\begin{equation}\label{connection}
{\textstyle
\gamma=\frac{(2A^2Ax-4A^2By+6ABAy+9BBx)}{4A^3+27B^2}dx+\frac{(4A^2Ay+6ABx+18BBy-9BAx)}{4A^3+27B^2}dy.}
\end{equation}
Note that, in general, the connection form becomes singular on the discriminant curve of the web. A direct calculation shows that the condition (\ref{associativity connection}) is equivalent to the fact that the connection  (\ref{connection}) remains holomorphic in the singular points. Therefore the germ of the connection form  is exact: $\gamma =df$, where $f$ is some function germ.

\begin{theorem}\cite{Ac}\label{homoexact}
Suppose ODE (\ref{binary}) admits an infinitesimal symmetry $X$ vanishing at the point $(0,0)$ on the discriminant curve $\Delta$ and the germ of the Chern connection form is exact  $\gamma=d(f)$,
where  $f$ is some function germ. Then  the equation germ and the symmetry are biholomorphic to one of the
following normal forms:
$$
\begin{array}{lll}
1) & y^{m_0}p^3-p=0, & X=(2+m_0)x\partial _x+2y\partial _y,\\
2) & p^3+2xp+y=0,    & X=2x\partial _x+3y\partial _y,\\
3) & (p-\frac{2}{3}x)(p^2+\frac{2}{3}xp+y-\frac{2}{9}x^2)=0, & X=x\partial _x+2y\partial _y,\\
4) & p^3+4x(y-\frac{4}{9}x^3)p+y^2+\frac{64}{81}x^6-\frac{32}{9}yx^3=0, & X=x\partial _x+3y\partial _y,\\
5) & p^3+y^2p=\frac{2}{\sqrt{27}}y^3\tan(2\sqrt{3}x), & X=y\partial _y, \\
6) & p^3+y^{3+m_0}p=y^{\frac{9+3m_0}{2}}F\left(\left[(m_0+1)\right]xy^{\frac{1+m_0}{2}}\right),& X=(1+m_0)x\partial _x-2y\partial _y,
\end{array}
$$
where $m_0$ is a non-negative integer and  $F(t)$ solves $$[12+2t^2-9tF]\frac{dF}{dt}=\frac{2(m_0+3)}{m_0+1}(4+27F^2)$$ with $F(0)=0$.
The weights $[w_1:w_2]$  uniquely determine  the normal form.
\end{theorem}
Characteristic 3-webs of associativity equations have an infinitesimal symmetry at each singular point and  holomorphic Chern connections in the normalization (\ref{connection}), therefore their singularities are equivalent to  the forms in Theorem \ref{homoexact}. Somewhere else we will address he question if each of the above singularities type realizes via characteristic 3-webs of associativity equations.

\section{Concluding Remarks}

\subsection{Generalization}
The geometric construction, presented in this paper, can be generalized to higher dimensions. As a result we obtain, for n-dimensional Frobenius manifold, a collection of $n$ commuting vector fields $v_i$ in $(\mathbb C^{n-1},0)$ satisfying the equation $\sum^n_{i=1}v_i=0$, i.e., a flat n-web germ of curves in $(\mathbb C^{n-1},0)$ admitting a "linear" symmetry.

\subsection{Acknowledgement}
 This research was partially supported  by the National Institute of Science
and Technology of Mathematics INCT-Mat.

\end{document}